\documentclass[12pt]{amsart}

\newtheorem{theorem}{Theorem}[section]

\newtheorem{corollary}[theorem]{Corollary}

\newtheorem{lemma}[theorem]{Lemma}
\newtheorem{proposition}[theorem]{Proposition}

\theoremstyle{definition}

\newtheorem{remark}[theorem]{Remark}

\newtheorem{example}[theorem]{Example}

\theoremstyle{parrafo}

\numberwithin{equation}{theorem}

\begin{document}

\title[]{A stability version of H\"older's inequality}

\author{J. M. Aldaz}
\address{PERMANENT ADDRESS: Departamento de Matem\'aticas y Computaci\'on,
Universidad  de La Rioja, 26004 Logro\~no, La Rioja, Spain.}
\email{jesus.munarrizaldaz@dmc.unirioja.es}

\address{CURRENT ADDRESS: Departamento de Matem\'aticas,
Universidad  Aut\'onoma de Madrid, Cantoblanco 28049, Madrid, Spain.}
\email{jesus.munarriz@uam.es}


\thanks{The author was partially supported by Grant MTM2006-13000-C03-03 of the
D.G.I. of Spain}







\begin{abstract} We present a stability version of H\"older's inequality,
incorporating an extra term that measures the deviation from
 equality. Applications are given.
\end{abstract}


\maketitle


\section
{Introduction.}

\markboth{J. M. Aldaz}{H\"older's inequality}

In the field of geometric inequalities, the expression {\em Bonnesen type} is used after Bonnesen classical refinement of the isoperimetric
inequality (cf., for instance, \cite{Os1}, \cite{Os2}), where the deviation from the case of equality (the disk) is given in terms of the outer radius and the inradius of
a bounded convex body.
The term {\em stability type} inequality is also used in a
related way (cf. \cite{Gr}), meaning that if  the deviation
from equality is ``small", then the objects under consideration must be
``close" to the extremal object.

Here we explore the question of
what a Bonnesen or stability version of H\"older's inequality
should look like, as we move away from the equality case.
Since the functions $f$ and $g$ involved in  H\"older's inequality will usually belong to different spaces, before
they can be compared we need to map these functions, with controlled distortion,
into a ``common measuring ground". The way we choose to do this is
by first normalizing, and then applying
 the Mazur map from $L^p$ and $L^q$ to $L^2$. For nonnegative functions
 in the unit sphere of $L^p$ the Mazur map
into $L^2$ is defined by $f\mapsto f^{p/2}$. We will be able to utilize
 its well known properties (cf. for instance,
 \cite{BeLi}) to obtain useful estimates.

 As a model for the stability version of H\"older's inequality, we use the (real) Hilbert space parallelogram
identity, suitably rearranged under the assumption that the
vectors are nonzero (see (\ref{bonhilb}) below).
With (\ref{bonhilb}) in mind we obtain a natural, straightforward
generalization of the parallelogram identity, valid for $1 < p < \infty$, though when $p\ne 2$ equality will of course be lost, cf. (\ref{bonhold}).
After one has decided which inequality to prove, the argument is
standard. In fact, it is {\em the} standard argument: From a refined
Young's inequality one obtains a refined H\"older inequality, which
in turn entails a refined triangle inequality,  which (together
with a simple additional observation) yields the uniform convexity
of $L^p$ spaces in the real valued
case, with optimal power type estimates for the modulus of convexity.

Like the parallelogram identity in the Hilbert space setting, (\ref{bonhold}) brings to the fore
the geometry of $L^p$ spaces, and conveys essentially the same information:
In order for $\|fg\|_1$ to be close to
$\|f\|_p \|g\|_q$, the angle between the
nonnegative $L^2$ functions
$|f|^{p/2}$ and $|g|^{q/2}$ must be small, with equality in
$\|fg\|_1\le \|f\|_p \|g\|_q$ precisely when the angle is zero.
Since
H\"older's inequality is one of the most often used inequalities, the refinement given here
is likely to have repercussions far beyond the few applications presented
below.

 The   paper is organized as follows. Section 2 contains the basic
inequality and its  proof, together with the
precedents I have been able to find, and a small discussion as to why
some plausible improvements of (\ref{bonhold}) cannot hold. Section 3
establishes a few direct consequences regarding bounds on interpolated norms. Specializing the previous remark about angles to the function 1 on
a probability space, we obtain a stability version of the following standard application of H\"older's inequality:
If $0 < r < s$, then
 every $f\in L^s$ satisfies $\|f\|_r\le \|f\|_s$, with equality if and
 only if  $|f|$ is constant. As we noted, the norms $\|f\|_s$ and  $\|f\|_r$ will be
close if and only if the angle between 1
and $|f|^{s/2}$ is small (cf.  Theorem \ref{interp1}).
Expressing this result in terms of the variance of $|f|^{s/2}$,
we shall see that
$\|f\|_s$ and  $\|f\|_r$ are close if and only if the normalized variance
 $\operatorname{Var} \left(|f|^{s/2}/\||f|^{s/2}\|_2\right)$
 is sufficiently small, cf.
  Corollary \ref{variance}. These results provide qualitative information
about the behavior of $L^p$ norms, which apparently had not been
 noticed before.
  Finally, Section 4 contains
 a sharpened triangle inequality, leading to the proof of uniform convexity
 announced
above.

We work on an arbitrary measure space $(X, \mathcal{A}, \mu)$, whose mention
will usually be omitted; to avoid trivialities
we assume that $\mu$ is not identically zero, and (when
dealing with uniform convexity) that $X$ contains
at least two points.

\section
{The basic inequality.}

In this paper $p$ and $q$ always denote conjugate exponents, i. e.,
$q=p/(p-1)$, and unless otherwise stated, it is understood that $f\in L^p$, $g\in L^q$ and neither function is zero almost everywhere.
To motivate the variant of H\"older's inequality given below, let us consider first the situation
in a real Hilbert space setting. From the parallelogram identity
\begin{equation}\label{par}
\|x + y\|^2 + \|x-y\|^2 = 2 \|x\|^2 + 2\|y\|^2
\end{equation}
we get, after expanding $\|x + y\|^2$, replacing $x$ by $tx$,
 taking $t=\|y\|/\|x\|$, and factoring $\|x\|\|y\|$, the equality
\begin{equation}\label{bonhilb} (x, y)  =
\|x\|\|y\|\left(1 - \frac12
\left\|\frac{x}{\|x\|}-\frac{y}{\|y\|}\right\|^2\right),
\end{equation}
valid for nonzero $x$ and $y$.
We follow  this line of thought in the $L^p$ setting, using
(\ref{bonhilb}) as a model. Observe that the identity
(\ref{bonhilb}) can be regarded as a stability version
(and also a proof) of the Cauchy-Schwarz inequality.

The first step is to refine Young's inequality $u^p/p + v^q/q - uv\ge 0$.

\begin{lemma}\label{betteryoungl}  Let $1 < p \le 2$ and let $q$ be its conjugate
exponent. Then for all $u,v \ge 0$
\begin{equation}\label{betteryoung}
\frac1q \left(u^{p/2} - v^{q/2}\right)^2  \le
\frac{u^p}{p} + \frac{v^q}{q} - uv   \le \frac1p \left(u^{p/2} - v^{q/2}\right)^2.
\end{equation}
\end{lemma}

\begin{proof} If $p=2 = q$ the result is trivial, so assume
 $1 < p < 2$. We prove the first inequality; the second can be obtained via an
essentially identical argument, by interchanging the roles of $p$ and $q$, and of $u$ and $v$. If either $u = 0$ or $v = 0$, formula
(\ref{betteryoung}) is obviously true. Fix $p$, fix $u > 0$, and
suppose $v > 0$. Expanding the square and simplifying, we see that
it is enough to check the following inequality:
\begin{equation} f(v):= \frac{2-p}{p} u^p + \frac{2}{q} u^{p/2} v^{q/2} - uv
\ge 0.
\end{equation}
Now  $v= u^{p-1}$ is the unique solution of $f^\prime (v)=0$. Since $f^{\prime\prime} > 0$, $f(u^{p-1}) =0$ is the global minimum  of  $f$.
\end{proof}

An extension of (\ref{bonhilb}) to the case $1 < p < \infty$ follows now by repeating the steps
in the usual derivation of H\"older's inequality from Young's inequality.
Only minimal modifications to the Hilbert space argument
given above
are needed, though of course, the equality becomes a two sided
 inequality when $p\ne 2$. We write $t_+ := \max\{t,0\}$ for the
 positive part of a real number or a real valued function, and
$t_+^r := \left(\max\{t,0\}\right)^r$, so the maximum is taken
first. The left hand side of the identity $\|f\|_p^{p/2} =
\||f|^{p/2}\|_2$ seems to be typographically more convenient and
easier to read than the right hand side, so  we will use it below.
However, it makes it less obvious that in (\ref{bonhold}) the
functions $|f|^{p/2}/\|f\|_p^{p/2}$ and $|g|^{q/2}/\|g\|_q^{q/2}$
are simply norm 1 vectors in $L^2$ (so we are in fact dealing with
the angle between $|f|^{p/2}$ and $|g|^{q/2}$, cf. Remark
\ref{angle}).

\begin{theorem}\label{betterhold}  Let $1 < p < \infty$ and let $q = p/(p-1)$ be its
conjugate exponent. If $f\in L^p $, $g\in L^q$, $\|f\|_p, \|g\|_q > 0$, and $ 1 < p \le 2$, then
\begin{equation}\label{bonhold}
\|f\|_p\|g\|_q \left(1 - \frac1p
\left\|\frac{|f|^{p/2}}{\|f\|_p^{p/2}}-\frac{|g|^{q/2}}{\|g\|_q^{q/2}}\right\|_2^2\right)_+\le\|fg\|_1 \le  \|f\|_p\|g\|_q \left(1 - \frac1q
\left\|\frac{|f|^{p/2}}{\|f\|_p^{p/2}}-\frac{|g|^{q/2}}{\|g\|_q^{q/2}}\right\|_2^2\right),
\end{equation}
while if $ 2 \le p < \infty$, the terms $1/p$ and $1/q$ exchange their positions in the preceding inequalities.
\end{theorem}
\begin{proof} Suppose $1 < p \le 2$. Write $u = |f(x)|$ and $v = |g(x)|$ in (\ref{betteryoung}), integrate, substitute $tf$ for $f$, and
set $t=\|g\|_q^{1/(p-1)}/\|f\|_p$. Now (\ref{bonhold}) immediately follows. If $2 \le p <\infty$, just interchange the roles of $p$
 and $q$.
\end{proof}

Of course,
when $p=2$ the inequality (\ref{bonhold}) follows from (\ref{bonhilb}), and
in fact, it is identical to it, save for the fact that only nonnegative
functions appear in (\ref{bonhold}).

The reason why we  take the positive part in the left hand side of (\ref{bonhold}), is that in some inequalities given
below we will need to take powers of the corresponding quantities.

\begin{remark}\label{angle} Recall that in a real inner product space, the angle $\angle (x,y)$ between  $x$ and $y$ is defined by
\begin{equation*}
\angle (x,y) : = \arccos \left(\frac{(x,y)}{\|x\|\|y\|}\right) =\arccos \left(1 - \frac12
\left\|\frac{x}{\|x\|}-\frac{y}{\|y\|}\right\|^2\right),
\end{equation*}
where the second equality follows from (\ref{bonhilb}). Actually,
the simpler expression $\theta (x,y) :=
\left\|\frac{x}{\|x\|}-\frac{y}{\|y\|}\right\|$, giving the length
of the segment between $x/\|x\|$ and $y/\|y\|$, is often taken as
the {\em definition} of angle in a general Banach space (cf., for
instance, pg. 403 of \cite{Cl}). 
In the real Hilbert space setting,
$\angle (x,y)$ and $\theta (x,y)$ are clearly comparable quantities
(in fact, $\theta (x,y) \le \angle (x,y) \le (\pi/2) \theta (x,y)$)
so up to a constant it does not matter which one is used. Thus, the
geometric content of (\ref{bonhold}) is clear: $\|fg\|_1 \approx
\|f\|_p \|g\|_q$ if and only if the angle $\angle
(|f|^{p/2},|g|^{q/2})$ is small. Note also that the same term
$\theta^2 (|f|^{p/2},|g|^{q/2})$ appears both on the left and on the
right hand sides of (\ref{bonhold}); hence, the exponent 2 cannot be
improved. This helps to explain why from (\ref{bonhold}) we obtain
optimal asymptotic power type estimates for the modulus of convexity
of $L^p (X, \mathbb{R})$ spaces.
\end{remark}

Observe that if $f$ and $g$ have disjoint supports then (\ref{bonhold})
becomes
\begin{equation}\label{half}
\|f\|_p\|g\|_q \left(1 - \frac2p\right)_+\le\|fg\|_1 = 0 \le  \|f\|_p\|g\|_q \left(1 - \frac2q\right).
\end{equation}
Hence, the right hand side bound worsens as $p\to 1$ (and $q\to\infty$).
Note also that the constant $1/2$  appears, instead of $1/p$ and
$1/q$, both in (\ref{bonhilb}) above and in
(\ref{pecaric}) below.
Thus, it is  natural to ask whether it is possible to
improve at least one of
the factors $1/p$, $1/q$ in (\ref{bonhold}), replacing it by
$1/2$ (of course, when supports are disjoint we cannot do better
than writing $0$ on the left hand side, but under less than full orthogonality, the change from $1/p$ to
$1/2$ might be useful).
  Next we show that such change is not possible.

\begin{example} Let $1 < p < 2$. Replacing $1/q$ by $1/2$ in the right hand side
of (\ref{bonhold}) and simplifying we find that this modification
of the second inequality
is equivalent to
\begin{equation}\label{nogo1}
\int|fg| \le   \|f\|_p^{1 - p/2} \|g\|_q^{1 - q/2} \int |f|^{p/2} |g|^{q/2}.
\end{equation}
Likewise, replacing $1/p$ by $1/2$ in the left hand side
of (\ref{bonhold}) leads to
\begin{equation}\label{nogo2}
\int|fg| \ge   \|f\|_p^{1 - p/2} \|g\|_q^{1 - q/2}
\int |f|^{p/2} |g|^{q/2}.
\end{equation}
It is easy to find examples showing that neither (\ref{nogo1})
nor (\ref{nogo2}) hold. Take for instance $f\equiv 1$ on $[0,1]$
and $g = 2 \chi_{[0,1/2]}$. Then $\|g\|_q = 2^{1 - 1/q}$, so
$1 = \int fg <  \|g\|_q^{1 - q/2}
\int  |g|^{q/2} = 2^{1/2 - 1/q}$ and thus (\ref{nogo2})
fails. Choosing now  $g\equiv 1$
and $f = 2 \chi_{[0,1/2]}$ we have $\|f\|_p = 2^{1 - 1/p}$ and
$1  > \|f\|_p^{1 - p/2}
\int |f|^{p/2}  = 2^{1/2 - 1/p}$ so (\ref{nogo1}) does not hold
either.

A more indirect argument shows that in fact $1/q$ cannot be
replaced by any fixed constant $c \in (0, 1/2)$ (independent of
$p$, or equivalently, of $q$). Since (\ref{bonhold}) can be used
to prove the uniform convexity of $L^p$ for $p >1$, if there were
such a $c$, then the
upper bound in (\ref{bonhold}) would not degenerate as $p\downarrow 1$,
and we would be able to show that the modulus of convexity of $L^p$ is
independent of $p$ for every $p\in (1,2]$, an obviously false result.
\end{example}

Despite its obvious interest, not much work has been done,
as far as I know,
regarding stability versions of H\"older's inequality. I am aware of two previous articles giving
bounds for the deviation from the case of equality. In
\cite{DraGo} the following result is presented:
\begin{equation}\label{drago}
0 \le 1 - \frac{(|f|,|g|)}{\|f\|_p\|g\|_q} \le \left(\frac{|f|^{p}}{\|f\|_p^{p}}-\frac{|g|^{q}}{\|g\|_q^{q}}, \frac1q \log |f| - \frac1p \log|g| \right)
\end{equation}
\begin{equation*}
\le\log\left[\frac{\left(|f|^{\frac{1 + pq}{q}}, |g|^{-\frac{1}{p}}\right) \left(|g|^{\frac{1 + pq}{p}}, |f|^{-\frac{1}{q}}\right)}{\|f\|_p^p\|g\|_q^q}\right],
\end{equation*}
where $(f,g):= \int fg$. Note that (\ref{drago}) does not coincide
with
the rearranged parallelogram identity (\ref{bonhilb}) when $p=q=2$.

An inequality more closely related to (\ref{bonhold}), which for nonnegative
 functions does extend (\ref{bonhilb}),
appears in \cite{PeSi}. The argument is actually the same as the
one used here (and in the standard proof of H\"older's inequality), save for the fact that the initial refinement of
Young's inequality is different from (\ref{betteryoung}). Suppose
$f, g\ge 0$. By Theorem
2 of \cite{PeSi}, if $1 < q \le 2 \le p <\infty$, then
\begin{equation}\label{pecaric}
\frac12
\frac{\left\|g^{2-q} \left( f\|g\|_q^{q/p} - g^{q-1} \|f\|_p \right)^2 \right\|_1}{\|f\|_p\|g\|_q^{q/p}}
\le \|f\|_p\|g\|_q  - \|fg\|_1 \le \frac12
\frac{\left\|f^{2-p} \left( g\|f\|_p^{p/q} - f^{p-1} \|g\|_q \right)^2 \right\|_1}{\|f\|_p^{p/q}\|g\|_q}.
\end{equation}
In addition to the factor $1/2$ mentioned before, there are other differences between  (\ref{pecaric}) and
(\ref{bonhold}). Note, for instance,
that every term in (\ref{bonhold}) is finite, while for $p > 2$,
whenever the support of $g$ is not contained in the support of $f$ the right hand side of (\ref{pecaric}) blows
up.

After submmiting this paper I have come accross the article
\cite{GGS}, where a refinement of H\"older's inequality is
obtained by using the positive definiteness of the Gram matrix. Write $m
:= \min \{ p^{-1}, q^{-1}\}$. Under the usual hypotheses,
 Theorem
2.3 of \cite{GGS} states that
\begin{equation}\label{gao}
\left( f, g\right) \le \|f\|_p\|g\|_q \left(1-r \right)^m,
\end{equation}
where $r$ is an explicitly defined function of $f^{p/2}, g^{q/2}$
and a third normalized vector $h\in L^2$. Both inequalities
(\ref{gao}) and (\ref{bonhold}) have in common the use of
$L^2$ to bound the deviation from equality. As differences, we
note that (\ref{gao}) is one sided, and it does not reduce to the
rearranged parallelogram identity when $p = q =2$.

Another relevant reference was found too late to include it in the
accepted version of the manuscript, cf. \cite{Si}. The one sided refinement given there is less related to (\ref{bonhold}) than those from \cite{DraGo},
\cite{PeSi}, and \cite{GGS}.

\begin{remark} It is easy to give a stability version
of the following standard variant of H\"older's inequality:
If $r > 0$,  $p^{-1} + q^{-1} = r^{-1}$, $f\in L^p$, and $g\in L^q$, then
$\|f g\|_r \le \|f\|_p \|g\|_q$. From it and an induction argument, stability versions
for multiple products can be obtained, that is, for the inequality
$\|\Pi_{i=1}^n f_i\|_r \le \Pi_{i=1}^n \|f_i\|_{p_i}$, where
$f_i\in L^{p_i}$ and $\sum_{i=1}^n p_i^{-1} = r^{-1}$.
\end{remark}

\section
{Interpolation-type consequences.}

In this section we derive some immediate interpolation-type results.
Note that
\begin{equation}
\left(1 - \frac1q
\left\|\frac{|f|^{p/2}}{\|f\|_p^{p/2}}-\frac{|g|^{q/2}}{\|g\|_q^{q/2}}\right\|_2^2\right) = 1 - \frac2q
\left(1 - \frac{\int |f|^{p/2} |g|^{q/2} }{\left(\int |f|^{p}\right)^{1/2}
\left(\int |g|^{q}\right)^{1/2}}\right),
\end{equation}
and these quantities are strictly positive when $q > 2$.
In what follows, both expressions will be used.

Recall  that
on a probability space, if $0 < r < s$, then
 every $f\in L^s$ satisfies $\|f\|_r\le \|f\|_s$, a fact that follows
 either from Jensen's inequality, or by writing $|f|$ as the product $|f|\cdot 1$ and then applying  H\"older's inequality. From the equality case
in either Jensen or H\"older inequalities, we have $\|f\|_r < \|f\|_s$ unless $|f|$ is
constant. This suggests that  the deviation of $|f|$ (or more precisely, of its normalized image under the Mazur map) from its mean
value can be used to obtain finer bounds.

\begin{theorem}\label{interp1}  Let $0 < r < s < \infty$, and let $f\in L^s$
satisfy
 $\|f\|_s > 0$. If $s \le 2r$, then
\begin{equation}\label{refint}
 \|f\|_s  \left[1 - \frac{2r}{s}
\left(1 - \frac{\||f|^{s/2}\|_1}{\||f|^{s/2}\|_2}\right) \right]_+^{1/r}
 \le \|f\|_r
 \le  \|f\|_s \left[1 - \frac{2(s-r)}{s}
\left(1 - \frac{\||f|^{s/2}\|_1}{\||f|^{s/2}\|_2}\right) \right]^{1/r},
\end{equation}
while if $s \ge 2r$, the  inequalities hold  with $2r/s$ and
$2(s - r)/s$ interchanged.
\end{theorem}
\begin{proof} We use Theorem (\ref{betterhold}) with $p = s/r > 1$,
$|f|^{r}\in L^p$, $q = s/(s-r) > 1$ and $g \equiv 1$.
Suppose first that $s \le 2r$, i.e., that $1 < p \le 2$. Substituting
in (\ref{bonhold}) and simplifying we get (\ref{refint}).
If $2\le p < \infty$ argue in the same way and use
the last part of Theorem \ref{betterhold}.
\end{proof}
A more common measure of the dispersion of $|f|^{s/2}/\||f|^{s/2}\|_2$
around its mean is the variance Var. From the previous result it is possible to derive bounds
for $\|f\|_r $ in
terms of
$\operatorname{Var}\left(|f|^{s/2}/\||f|^{s/2}\|_2\right)$.

\begin{corollary}\label{variance}  Let $0 < r < s < \infty$, and suppose
 $0 < \|f\|_s < \infty$. If $s \le 2r$, then
\begin{equation}\label{int1var}
\|f\|_s  \left[1 - \frac{2r}{s}
\operatorname{Var} \left(\frac{|f|^{s/2}}{\||f|^{s/2}\|_2}\right) \right]_+^{1/r}
 \le \|f\|_r
 \le  \|f\|_s  \left[1 - \frac{s - r}{s}
\operatorname{Var}\left(\frac{|f|^{s/2}}{\||f|^{s/2}\|_2}\right) \right]^{1/r},
 \end{equation}
while if $s \ge 2r$, the  same inequalities hold, but
 with the terms $2r/s$ and $(s-r)/s$ interchanged.
 \end{corollary}

\begin{proof} Note that for all $x\in [0,1]$
\begin{equation}\label{eleminequ}
 2^{-1} (1 - x^2) = 2^{-1} (1 + x) (1-x) \le  1 - x \le   1 - x^2.
\end{equation}
Next we set $x = \||f|^{s/2}\|_1/\||f|^{s/2}\|_2$. Then
$x \le 1$ by either Jensen's inequality or more simply, the
nonegativity of the variance. Substituting in (\ref{eleminequ}) we
obtain
\begin{equation}\label{var}
\frac{1}{2}
\operatorname{Var}\left(\frac{|f|^{s/2}}{\||f|^{s/2}\|_2}\right) \le
 \frac{\||f|^{s/2}\|_2 - \||f|^{s/2}\|_1}{\||f|^{s/2}\|_2}
  \le
\operatorname{Var} \left(\frac{|f|^{s/2}}{\||f|^{s/2}\|_2}\right),
 \end{equation}
Now (\ref{int1var}) follows from (\ref{refint}) when $s/r \le 2$,
while if if $2\le s/r $, we use the last part of Theorem \ref{interp1}
to obtain the corresponding inequalities.
\end{proof}

Theorem \ref{interp1} and Corollary \ref{variance} are stability
results, in the sense that   $\|f\|_s$ and  $\|f\|_r$ are ``close"
if and only if $|f|$ is ``nearly" constant;  when $\angle
(|f|^{s/2}, 1)$ (or $\operatorname{Var}
\left(|f|^{s/2}/\||f|^{s/2}\|_2\right)$) is sufficiently small,
these norms are comparable. We believe these results will be useful
in contexts where information
 is available about
the first and second moments of a function, as is often
the case in Probability Theory.

\begin{remark} It is easy to check that the factors between square brackets in the left hand sides
of (\ref{refint}) and (\ref{int1var}) can actually be negative, so
the positive part must be taken before raising them to the $1/r$ power.
Take for instance, $s =2$, any fixed $r \in (1,2)$,  and $f = \sqrt{n} \chi_{[0,1/n]}$ on $[0,1]$, with $n = n(r)$ ``large enough".
\end{remark}

A variant of the result on containment of $L^p$ spaces
exchanges the probability measure (or more generally, finite measure)
hypothesis by the condition that $f$ belongs to $L^{p_0}$, for some
$p_0 < p$. We consider this next.

\begin{theorem}\label{interp2}  Let $0 < p_0 < p < p_1 < \infty$, and let
$t = t(p)$ be given by the equation $p^{-1} = (1 -t) p_0^{-1} + t p_1^{-1}$. Suppose
$f\in L^{p_0} \cap L^{p_1}$ and $f\ \slash\hspace{-0.38cm}\equiv 0$.
If $p_0/p_1 \le t^{-1} - 1$, then
\begin{equation}\label{formulaint2}
 \|f\|_{p_0}^{1-t}  \|f\|_{p_1}^{t} \left[1 - \frac{2(1-t) p_1}{(1-t) p_1 + t p_0}
\left(1 - \frac{\int |f|^{\frac{p_0 + p_1}{2}}}{\left(\int|f|^{p_0}\right)^{1/2} \left(\int|f|^{p_1}\right)^{1/2}}\right) \right]_+^{1/p}
 \end{equation}
\begin{equation}\label{formulaint2a}
 \le \|f\|_p
 \le   \|f\|_{p_0}^{1-t}  \|f\|_{p_1}^{t} \left[1 -
 \frac{2tp_0}{(1-t) p_1 + t p_0} \left(1 - \frac{\int |f|^{\frac{p_0 + p_1}{2}}}{\left(\int|f|^{p_0}\right)^{1/2} \left(\int|f|^{p_1}\right)^{1/2}}\right) \right]^{1/p},
\end{equation}
while if $p_0/p_1 \ge t^{-1} - 1$, the  inequalities are reversed, and
 the positive part of the term between square brackets
 is taken in the right hand side of (\ref{formulaint2a}).
\end{theorem}
\begin{proof} Again we use Theorem (\ref{betterhold}),  with the
functions $f^{(1-t) p} f^{t p} = f^p$, and the conjugate
 exponents $p_0 /[(1-t)p]$ and $p_1/t p$.  Note that
 $p_0 /[(1-t)p] > 1$ and $ p_1/t p > 1$, while $p_0/p_1 \le t^{-1} - 1$ if and only
 if $p_0 /((1-t)p) \le 2$.
\end{proof}

\begin{remark} The preceding theorem leads to a midpoint interpolation result for arbitrary
pairs of functions. Suppose, for instance, that $f, h \in L^{p_0} \cap L^{p_1}$, $f, h \ \slash\hspace{-0.4695cm}\equiv 0$, $\|f\|_{p_0}\le \|h\|_{p_0}$, and
$\|f\|_{p_1}\le \|h\|_{p_1}$. It is easy to see that
$\|f\|_{p} > \|h\|_{p}$ may happen for some intermediate $p\in (p_0,p_1)$.
Consider the following example: Set $f(x) = (1 - 1/n)\chi_{[0,1/2]}$ on $[0,1]$, where $n \ge 6$ is fixed, and let
$h(x) = x$.
Then $\|f\|_{1} < \|h\|_{1}$ and $\|f\|_{\infty} < \|h\|_{\infty}$,
but $\|f\|_{n} > \|h\|_{n}$.
Note that $\|f\|_{p} < \|h\|_{p}$ for every large enough $p < \infty$; in particular, if $n = 6$ we can take $p_1 = 11$, so there is a reversal
of the inequality at $p=(p_0 + p_1) /2$.
 However,  under the additional condition
on the angles
$\angle (|h|^{p_0/2}, |h|^{p_1/2}) \le \angle (|f|^{p_0/2}, |f|^{p_1/2})$,
or equivalently,  $\theta (|h|^{p_0/2}, |h|^{p_1/2}) \le \theta (|f|^{p_0/2}, |f|^{p_1/2})$, at the midpoint $p = (p_0 + p_1)/2$
 we have $\|f\|_{p} \le \|h\|_{p}$
whenever $\|f\|_{p_0}\le \|h\|_{p_0}$ and
$\|f\|_{p_1}\le \|h\|_{p_1}$. To see this, note
that if $p = (p_0 + p_1)/2$, then $t = p_1/(p_0 + p_1)$, so from
(\ref{formulaint2a}) and (\ref{formulaint2}) we get
\begin{equation*}  \|f\|_p
 \le   \|f\|_{p_0}^{1-t}  \|f\|_{p_1}^{t} \left[1 -
 \frac{1}{2} \left\|\frac{|f|^{p_0/2}}{\left\||f|^{p_0/2}\right\|_{2}}  - \frac{|f|^{p_1/2}}{\left\||f|^{p_1/2}\right\|_{2}}\right\|_2^{2}\right]^{1/p}
\end{equation*}
 \begin{equation*}
 \le  \|h\|_{p_0}^{1-t}  \|h\|_{p_1}^{t} \left[1 -
 \frac{1}{2} \left\|\frac{|h|^{p_0/2}}{\left\||h|^{p_0/2}\right\|_{2}}  - \frac{|h|^{p_1/2}}{\left\||h|^{p_1/2}\right\|_{2}}\right\|_2^{2} \right]^{1/p}
 \le \|h\|_p.
 \end{equation*}
Needless to say, stronger assumptions on the angles lead to stronger
interpolation results. For instance, if
$\theta (|h|^{p_0/2}, |h|^{p_1/2}) < \theta (|f|^{p_0/2}, |f|^{p_1/2})$,
then  $\|f\|_{p} < \|h\|_{p}$ for every $p$ in some neighborhood of
$(p_0 + p_1)/2$, since the quantities involved in (\ref{formulaint2a}) and (\ref{formulaint2}) change continuously. It is also possible to consider
conditions of the type $\|f\|_{p_i}\le c_i\|h\|_{p_i}$, with $c_i > 0$ not necessarily  equal to 1, or even
to have  $h \in L^{r_0}\cap L^{r_1}$ with $r_i \ne p_i$, as is often
done in interpolation theorems. But we will not pursue these elaborations
here.
\end{remark}

\begin{remark} In standard interpolation results, such as the Riesz-Thorin and the Marcinkie\-wicz interpolation theorems, the pairing
between the functions $f$ and $h = T(f)$
is not arbitrary but given respectively by a linear or sublinear operator $T$, and the conclusion, of
course, is much stronger than anything contained in the previous
remark. The attentive reader may wonder why more general pairings are interesting, or in other words, whether there
is any need to go beyond sublinearity. Next we give an example
where such a result might be useful. It involves
the derivative $DMf$ of the one dimensional, uncentered Hardy-Littlewood maximal
function $Mf$, defined as follows:
 Given  a locally integrable function
$f:\mathbb{R}\to \mathbb{R}$,
$$
Mf(x) := \sup_{ x\in I}\frac{1}{|I|}\int_I |f(y)|dy,
$$
where $I$ is any interval containing $x$ and $|I|$ stands for its length.
Starting with the paper \cite{Ki}, there has been in recent years a growing
interest regarding the regularity of the maximal function (cf., for instance,
\cite{AlPe} and the references contained therein). Suppose for simplicity
that $f:\mathbb{R}\to \mathbb{R}$ is a  compactly supported
Lipschitz function. It is shown in \cite{Ki} (cf. also \cite{HaOn}) that
for every $1 < p \le \infty$ there is a constant $c_p$ (independent of
$f$) such that
$\|DMf\|_p \le c_p \|Df\|_p$. However, the methods used in
\cite{Ki} and \cite{HaOn} cannot tell us whether we actually have
$c_p < 1$, that is, whether the maximal operator
$M$ has a smoothing effect on $f$. For $p=1$, Theorem 2.5 of
\cite{AlPe} states  that $\|DMf\|_1 \le  \|Df\|_1$,
and $c_1 = 1$ is sharp, while for $p=\infty$, we have
 $\|DMf\|_\infty \le (\sqrt{2} -1) \|Df\|_\infty$ and $c_\infty =
(\sqrt{2} -1)$ is best possible, by
 \cite{ACP}. Thus, it is natural to conjecture ``by interpolation" that
whenever $1 < p <\infty$, the optimal constant $c_p$ satisfies
$c_p < 1$, and furthermore,
$\lim_{p\to\infty} c_p = \sqrt{2} -1$. Nevertheless, since the operator
$Df\mapsto DMf$ is neither linear nor sublinear, it falls outside
the realm of currently available interpolation theorems. Unfortunately,
the second endpoint for which information is available
happens to be $p = \infty$, so our stability version of
H\"older's inequality also fails to yield anything new on this
question.
\end{remark}

\section
{The triangle inequality and uniform convexity.}

 Like Clarkson's inequalities and Hanner's inequalities, formula
(\ref{bonhold}) can lay claim to being an $L^p$ generalization of
the parallelogram identity. Furthermore, despite its easy proof, the
refinement of H\"older's inequality presented above does have
strength: It gives, by sharpening Minkowski's inequality, the uniform
convexity of $L^p$ spaces (at least in the real valued case), with  the right asymptotic behavior of
the modulus of convexity for all  $p\in (1,\infty)$. The exact
asymptotic behavior was found by O. Hanner (cf. \cite{Ha}, or
\cite{LiTza2}, p. 63);  Clarkson's original inequalities (see the
Corollary in pg. 403 of \cite{Cl})  yield it  over the range $2\le p
< \infty$, but not for $1 < p < 2$.

The
arguments presented here only cover the real valued case, and the
complex valued case if $p\ge 2$. Since only the moduli of functions
(and not their signs) play any role  in the sizes of $\|fg\|_1$ and
$\|f\|_p \|g\|_q$, the same must necessarily happen with the error
terms in any refinement of H\"older's inequality. In particular,
this is the case with (\ref{bonhold}). But for some applications,
such as a refined triangle inequality, it would be preferable to
control the departure from maximal size in terms of $|f-g|$  rather
than  $||f| - |g||$. We shall show that for real valued functions,
and for complex valued functions when $p \ge 2$, one can assume the
comparability of $\|f-g\|_p$  and  $\|f| - |g|\|_p$. But the proof
in the complex case when $1 < p < 2$ has eluded us. A recent, new
proof of uniform convexity, relying on the notion of thin slices and
which does apply to the complex case, can be found in \cite{HaO}
(however, there the author is unconcerned about
 the precise behavior of the modulus of convexity).

The improved
Minkowski's inequality given next is obtained from our refinement of H\"older's inequality
by the usual ``duality" argument. By the
``duality" argument we do not mean knowing that the dual of $L^p$ is $L^q$, but simply that
\begin{equation}\label{dualnorm}
\|f\|_p =\sup_{\{g\in L^q: \|g\|_q = 1\}} \int fg,
\end{equation}
which follows from H\"older's inequality
together with the trivial observation that equality
is achieved when $g = \frac{|f|^{p - 1} \overline{\operatorname{sign} f}}{\|f\|_p^{p - 1}}$. Here $\operatorname{sign} (z) := e^{i\theta}$ for
every complex nonzero $z = r e^{i\theta}$, and
$\operatorname{sign} (0) := 1$ (we adopt this convention, rather than the usual $\operatorname{sign} (0) := 0$, since in order to multiply quantities
without changing sizes it is useful to always have
$|\operatorname{sign} (z)| = 1$). As is well known, (\ref{dualnorm}) immediately entails the triangle inequality:
\begin{equation}\label{proof}
\|f + h\|_p = \sup_{\{g\in L^q: \|g\|_q = 1\}} \int (f + h) g \le \sup_{\{g_1\in L^q: \|g_1\|_q = 1\}} \int fg_1 + \sup_{\{g_2\in L^q: \|g_2\|_q = 1\}} \int h g_2 = \|f\|_p + \|h\|_p.
\end{equation}
However, usually this proof appears with the explicit
maximizing $g$ written in place of the first supremum, and then it proceeds
from there. As it turns out, it
will be more convenient for us to do likewise below.

\begin{theorem}\label{triangle}  Let $1 < p < \infty$.  If $f, h\in L^p $, $\|f\|_p, \|h\|_p > 0$, and $ 1 < p \le 2$, then
\begin{equation}\label{bonmink}
\|f + h\|_p \le \|f\|_p \left(1 - \frac1q
\left\|\frac{|f + h|^{p/2}}{\|f + h\|_p^{p/2}}-\frac{|f|^{p/2}}{\|f\|_p^{p/2}}\right\|_2^2\right) + \|h\|_p \left(1 - \frac1q
\left\|\frac{|f + h|^{p/2}}{\|f + h\|_p^{p/2}}-\frac{|h|^{p/2}}{\|h\|_p^{p/2}}\right\|_2^2\right),
\end{equation}
while if $ 2 \le p < \infty$ the same inequality holds, but with
$1/p$ replacing $1/q$ throughout.
\end{theorem}
\begin{proof} Suppose $1 < p \le 2$. Then
\begin{equation*}
\|f + h\|_p = \int \frac{|f + h|^{p-1}}{\||f + h|^{p-1} \|_q} |f + h| \le
\int \frac{|f + h|^{p-1}}{\||f + h|^{p-1} \|_q} |f| +
\int \frac{|f + h|^{p-1}}{\||f + h|^{p-1} \|_q} |h|
\end{equation*}
and the result follows by applying (\ref{bonhold}).
If $2\le p < \infty$ argue in the same way and use
the last part of Theorem \ref{betterhold}.
\end{proof}

Next, we recall some basic facts about the Mazur map $\psi_{r,s} : L^r\to L^s$.
It is defined first on the
unit sphere by
$\psi_{r,s} (f) := |f|^{r/s} \mbox{ sign }f$, and then extended
to the rest of the space by homogeneity (cf. \cite{BeLi}, pp. 197--199 for additional
information on $\psi_{r,s}$). The ``angle"
$
\left\|\frac{|f|^{p/2}}{\|f\|_p^{p/2}}-\frac{|g|^{q/2}}{\|g\|_q^{q/2}}\right\|_2
$
in (\ref{bonhold})
is obtained by
 applying the Mazur maps from the nonnegative functions in the unit spheres of $L^p$ and $L^q$,
into the unit sphere of $L^2$. Thus, we have  control over the
distortion, since when $r < s$, the map $\psi_{s,r}$ is  Lipschitz
on the unit sphere of $L^s$, with
 constant $s/r$, while its inverse $\psi_{r,s}$ is H\"older
with exponent $r/s$. This is the content of the following well known lemma, included here for the reader's convenience.
It is a special case of Proposition 9.2, pp. 198-199
of \cite{BeLi}, cf. also the proof of Theorem 9.1, pg. 198, partially
sketched below. Note however that in \cite{BeLi} the harder, complex
valued case is handled, and the H\"older constant (as opposed to the
H\"older exponent) is not specified. We will consider
 the Mazur map acting only on nonnegative functions, since that  is all we shall use.
 In this easy case we show that the H\"older constant is 1.

\begin{lemma}\label{mazur}  Let $1 < r < s < \infty$, and let
$f,h \ge 0$. If $f, h\in L^r$ satisfy  $\|f\|_r = \|h\|_r = 1$, then  $\|f^{r/s} - h^{r/s}\|_s \le \|f - h\|_r^{r/s}$,
while if
$f, h\in L^s$ have norms  $\|f\|_s = \|h\|_s = 1$, then
$\|f^{s/r} - h^{s/r}\|_r \le (s/r) \|f - h\|_{s}$.
\end{lemma}
\begin{proof}  To prove the H\"older assertion, note that
by concavity of $t^\alpha$ for $0<\alpha < 1$, if $a > b$, then
$a^\alpha - b^\alpha \le (a-b)^\alpha$. Suppose $f$ and $h$ are nonnegative functions of norm 1 in $L^r$. Taking $\alpha = r/s$ and integrating the pointwise inequality
$|f^{r/s}(x) - h^{r/s}(x)|^s \le |f (x) - h (x)|^r$ we get
$\|f^{r/s} - h^{r/s}\|_s \le \|f - h\|_r^{r/s}$.

We sketch the proof the Lipschitz claim, directing the
reader to \cite{BeLi} for additional details. Let us denote by $d \psi_{s,r} (f) (h)$
the Gateaux (i.e., the directional) derivative of the Mazur map based at the point
$f$ and in the direction of $h$, where the nonnegative functions $f$ and $h$ belong  the unit sphere of $L^s$. It is enough to show that
$\|d \psi_{s,r} (f) (h)\|_r^r \le (s/r)^r$, which follows by explicit
computation of the directional derivative, and an application of H\"older's inequality together with $\|f\|_s = \|h\|_s =1$.
\end{proof}

After proving a simple lemma, we use the the properties of the Mazur map to express the preceding triangle inequality in terms of the $p$ norm.

\begin{lemma}  Let $x, y, z$ be vectors in a normed space,
and let $p\in (1,\infty)$.
Then $\|x - y\|^p \le 2^{p-1} \left(\|x - z\|^p +  \|y - z\|^p\right)$.
\end{lemma}
\begin{proof} We may assume that $x \ne y$. Since
$\|x - y\| \le  \|x - z\| + \|y - z\|$, writing $a:= \|x - z\|/\|x - y\|$ and $b:= \|y - z\|/\|x - y\|$ we have that $a + b \ge 1$ and
$(a ^p + b^p)\|x - y\|^p =  \|x - z\|^p + \|y - z\|^p$. Minimizing
$a ^p + b^p$ subject to  $a + b \ge 1$ and $a, b \ge 0$, we obtain $a = b =1/2$,
from which the result follows.
\end{proof}

Clarkson gave asymptotic estimates
for the modulus of convexity of $L^p$ of order
 $O(\varepsilon^p)$ when $2 \le  p < \infty$
and $O(\varepsilon^q)$ when $1 <  p \le 2$,  where
$\varepsilon = \|f - h\|_p$. The optimal estimate
$O(\varepsilon^2)$ when $1 <  p \le 2$ was found by Hanner.
It is easy for us  to explain this different behavior in terms of the
Mazur map: When $p \le 2$ the map $\psi_{2,p}$ is Lipschitz, and
hence the exponent $2$ in the error term from (\ref{bonhilb}) or
(\ref{bonhold}) is preserved, while if $p \ge 2$, then  $\psi_{2,p}$ is
$2/p$-H\"older, so the exponent 2 changes to $p$.

\begin{corollary} \label{trianpos} Let $1 < p < \infty$, and let $f, h\in L^p $. If
 $ 1 < p \le 2$, then
\begin{equation}\label{bonmink1}
\|f + h\|_p \le \|f\|_p + \|h\|_p  - \min \{\|f\|_p, \|h\|_p\}
 \left(\frac{p(p-1)}{8}
\left\|\frac{|f|}{\|f\|_p}-\frac{|h|}{\|h\|_p}\right\|_p^2\right),
\end{equation}
while if $ 2 \le p < \infty$,
\begin{equation}\label{bonmink2}
\|f + h\|_p \le \|f\|_p + \|h\|_p  - \min \{\|f\|_p, \|h\|_p\}
 \left(\frac{1}{2p}
\left\|\frac{|f|}{\|f\|_p}-\frac{|h|}{\|h\|_p}\right\|_p^p\right).
\end{equation}
\end{corollary}
\begin{proof} The result follows from
Theorem \ref{triangle}, the previous Lemma, and Lemma \ref{mazur}.
\end{proof}

Suppose, in order to simplify the corresponding expressions,
that $ \|f\|_p = \|h\|_p = 1$. A drawback of the preceding
corollary is that in the right hand side we have $\left\||f|-|h|\right\|_p$
rather than $\left\|f - h\right\|_p$, while the left hand side
depends on $f + h$, not on $|f| + |h|$. This is unavoidable since
we are deriving the result from the stability version of H\"older's inequality (\ref{bonhold}). Thus, the case where
$\left\||f|-|h|\right\|_p <<
\left\|f - h\right\|_p$ must be handled via a separate argument,
which somehow we have failed to find when $f$ and $h$ are complex valued
and $p < 2$. The real valued case is easy since the only possibility for
cancellation is to have opposite signs, and for $p\ge 2$
the complex valued case immediately
follows from the convexity of $p/2$.

 Note that the bound in the next proposition
  has nothing to do with uniform convexity:
It holds even when $p=1$. In fact, all we are doing is checking
the intuitively obvious fact that if we want $\|f + h\|_p$ to be large,
the signs of $f$ and $h$ must be very similar, specially if $p$ is small.
 While this ought to be also true in the complex valued
case, as I said I have not been able to prove it.

\begin{proposition}\label{noncomparison}  Let $1 \le p < \infty$, let $0 < t < 1$, and let  $f, h\in L^p$
be real valued functions. If   $\||f| - |h|\|_p^p < t \|f -
h\|_p^p$, then $\|f + h\|_p < \left(\left(\|f\|_p + \|h\|_p\right)^p
- (1 -t) \|f - h\|_p^p\right)^{1/p}$.
\end{proposition}
\begin{proof} First, we may assume  that $f \ge 0$, since by the convention $\operatorname{sign} (0) := 1$ (adopted just after (\ref{dualnorm}))
given any $x$ we have
$|f(x) - h(x) | =|f(x) \mbox{ sign } f(x) - h(x)\mbox{ sign } f(x)|$, and likewise for $|f(x) + h(x)|$. Next, note that if $a \ge 0$ and
$b\in \mathbb{R}$, then $|a + b|^p + |a-b|^p = |a + |b||^p + |a-|b||^p$, so
writing $f(x) = a$, $h(x) = b$, and integrating, we get
\begin{equation}\label{compare}\|f + h\|_p^p =
\|f +|h|\|_p^p + \|f - |h|\|_p^p - \|f - h\|_p^p
 \le
\left(\|f\|_p + \|h\|_p\right)^p   - (1 -t)  \|f - h\|_p^p.
\end{equation}
\end{proof}

\begin{remark}  Note that by Taylor's formula (or by
linear approximation at 0 and concavity), we have
$(1 - x)^{1/p} \le 1 - p^{-1} x$. Applying this inequality to the conclusion
of the previous proposition when $\|f\|_p = \|h\|_p
= 1$, we get
\begin{equation}\label{taylor}\left\|\frac{f + h}{2}\right\|_p
 \le
1  - \frac{1 -t}{p 2^p}  \|f - h\|_p^p.
\end{equation}
\end{remark}

Let $B$ be  a Banach space.
Clarkson's original definition of uniform convexity
 requires that for every $0 < \varepsilon \le 2$ there exist a
$\delta (\varepsilon) > 0$ such that if $\|f\| = \|h\| = 1$
and $\|f -h\| \ge\varepsilon$, then $\left\|\frac{f+h}2\right\| \le 1 - \delta (\varepsilon)$ (c.f Definition 1., pp.396-397 of \cite{Cl}).
The often used and seemingly weaker assumption $\|f\|, \|h\|
\le 1$ is of course equivalent to $\|f\| = \|h\|
= 1$, since $f$ and $h$ must have norm one in order
to maximize $\|f + h \|_p$ subject to
$\|f - h\|\ge\varepsilon$ (see Lemma 5.1 pg. 381 of \cite{Da}  for a full proof).
In the words of \cite{BaCaLi}, $B$ is uniformly convex if its unit
ball is ``uniformly free of flat spots".
From the viewpoint of the geometry of $B$ is is often interesting
to have a good estimate of how $\delta$ depends on $\varepsilon \in (0,2]$. The following definitions and results
are taken from \cite{LiTza2}, specially pg. 63.
The modulus of convexity $\delta_B$ of $B$ is given by
\begin{equation}\label{defmc}
\delta_B (\varepsilon) :=
\inf\left\{1 - \left\|\frac{f+ h}2\right\|: \|f\|=\|h\|= 1, \left\|f-h\right\| = \varepsilon \right\}.
\end{equation}
We say that $\delta_B$ is of power type $r$ if there exists a
constant $c > 0$ such that $\delta_B (\varepsilon) \ge c
\varepsilon^r$. For $B = L^p$ and $ 1 < p \le 2$, $\delta_B
(\varepsilon) = (p-1) \varepsilon^2/8 + o (\varepsilon^2)$, while
for $2 \le p <\infty$, $\delta_B (\varepsilon) =  \varepsilon^p/(p
2^p) + o (\varepsilon^p)$.

The next result shows that in the real valued case, the preceding
variants of the triangle inequality yield  the optimal value of $r$
in the power type estimates. The constants, however, are not
optimal. But they are not too far away from optimality either. We
make an effort to obtain ``fairly good" constants for the modulus of
convexity (and not just good power type estimates, which is all one
usually needs for applications) since this entails that the
constants in the original inequality (\ref{bonhold}) must also be
``fairly good".

\begin{theorem}\label{clarkson}  Let $1 < p < \infty$. Then $B = L^p(X, \mathbb{R})$ is uniformly
convex. Furthermore, its modulus of convexity satisfies the
following inequalities. If $p\in (1,2]$, then for every $c
> 1$
 there exists an $\varepsilon = \varepsilon (c)$
such that for all  $f, h\in L^p$ with $\|f\|_p = \|h\|_p = 1$ and
$\left\|f-h\right\|_p \le \varepsilon$,
\begin{equation}\label{ucsmp}
\delta_B (\left\|f-h\right\|_p) \ge \frac{p(p-1)}{16 c}
\left\|f-h\right\|_p^2.
\end{equation}
On the other hand, if $2\le p < \infty$, then for all  $f, h\in L^p$
with $\|f\|_p = \|h\|_p = 1$,
\begin{equation}\label{ucbigp}
\delta_B (\left\|f-h\right\|_p) \ge \frac{\left\|f-h\right\|_p^p}{p
2^p + 4p}.
\end{equation}
\end{theorem}

\begin{proof} Note that for every $p\in (1, \infty)$ and every $t\in
(0,1)$,  if $\||f| - |h|\|^p_p < t \|f - h\|^p_p$, then
\begin{equation}\label{ucbigp2}
\delta_B (\left\|f-h\right\|_p) \ge \frac{(1 - t)
\left\|f-h\right\|_p^p}{p 2^p}
\end{equation}
by Proposition \ref{noncomparison}, or more precisely, by
(\ref{taylor}).

We prove (\ref{ucbigp}) first. Given $t\in (0, 1)$, if $\||f| -
|h|\|^p_p \ge t \|f - h\|^p_p$, by (\ref{bonmink2}) we have the
bound
\begin{equation}\label{ucbigp1}
\delta_B (\left\|f-h\right\|_p) \ge \frac{t
\left\|f-h\right\|_p^p}{4 p}.
\end{equation}
Choosing $t$ so that the lower bounds given by (\ref{ucbigp2}) and
(\ref{ucbigp1}) are equal, (\ref{ucbigp}) follows.

With respect to (\ref{ucsmp}), observe that for every $t\in (0,1)$
and $\left\|f-h\right\|_p$ sufficiently small (depending on $t$),
 the bound
\begin{equation}\label{ucsmp1}
\delta_B (\left\|f-h\right\|_p) \ge \frac{t^{2/p} p (p-1)
\left\|f-h\right\|_p^2}{16},
\end{equation}
which follows from  (\ref{bonmink1}) when $\||f| - |h|\|^p_p \ge t
\|f - h\|^p_p$, is always smaller than the bound given by
(\ref{ucbigp2}) when $\||f| - |h|\|^p_p < t \|f - h\|^p_p$. Writing
$c = t^{-2/p}$, (\ref{ucsmp}) follows by fixing $\varepsilon > 0$
small enough and taking $\left\|f-h\right\|_p \le \varepsilon$.
\end{proof}

We have given an asymptotic estimate when $1 < p \le 2$ in order to
be as precise as we can. If we are not concerned with good
constants, to obtain a statement which does not require
$\varepsilon$ to be small we can just fix any $t$ (say $t = 2^{-1}$
for definiteness) and take the minimum of the quantities given by
(\ref{ucbigp2}) and (\ref{ucsmp1}).

Next we consider the case of $L^p(X, \mathbb{C})$ spaces, when $p\ge
2$. The argument is essentially the same as in Proposition
\ref{noncomparison}.

\begin{proposition}\label{compnoncomparison}  Let $2 \le p < \infty$,
let $0 < t < 1$, and let  $f, h\in L^p$ be complex valued functions.
If $\||f| - |h|\|_p^p < t \|f - h\|_p^p$, then
$\|f + h\|_p <
\left(\left(\|f\|_p + \|h\|_p\right)^p  - (1-t) \|f -
h\|_p^p\right)^{1/p}$.
\end{proposition}
\begin{proof} As before, we may assume  that $f \ge 0$. Writing
$h = |h| e^{i\alpha}$, where $\alpha = \alpha (h(x))$, we have that  for every $x$,
\begin{equation} \label{lhs}
|f(x) + h(x)|^p + |f(x) - h(x) |^p = \end{equation}
\begin{equation*}
|f^2(x) + |h(x)|^2  + 2 f(x) |h(x) |\cos \alpha(h(x))|^{p/2} +
|f^2(x) + |h(x)|^2  - 2 f(x) |h(x) |\cos \alpha(h(x))|^{p/2}.
\end{equation*}
By the convexity of $t^{p/2}$,
\begin{equation} \label{largep}
|f(x) + h(x)|^p + |f(x) - h(x) |^p
\end{equation}
\begin{equation*} \le
|f^2(x) + |h(x)|^2  + 2  |h(x) |f(x)|^{p/2} +
|f^2(x) + |h(x)|^2  - 2 |h(x) |f(x)|^{p/2}
\end{equation*}
\begin{equation*} = |f(x) + |h|(x)|^p + |f(x) - |h|(x) |^p .
\end{equation*}
The rest of the proof is as in Proposition \ref{noncomparison}.
\end{proof}

\begin{remark} From the preceding Proposition and the second part of Corollary
\ref{trianpos},
 the uniform convexity of the $L^p(X, \mathbb{C})$ spaces
when $p\ge 2$ follows in exactly the same way and with the same
constants as in Theorem \ref{clarkson}, so we avoid the repetition.
\end{remark}

\begin{remark} As we have noted, a disadvantage of the refined triangle
inequality given in Corollary \ref{trianpos}, is that the error or
stability term depends only on the moduli of the functions involved,
and not their signs. But this inequality has its advantages also.
One of them is that it interacts well with other inequalities given
here, in the sense that it is easy to obtain nontrivial information
by combining them. For instance, suppose $\mu (X) = 1$ and $0 < r <
s$, with $f,h\in L^s$. Under suitable hypotheses on the variance of
$|f + h|^{s/2}$, we can easily find bounds for $\|f + h\|_{s}$ in
terms of $\|f\|_r$ and $\|h\|_r$, by using Theorem \ref{interp1} or
Corollary \ref{variance}, together with Corollary \ref{trianpos}.
Alternatively, we might be interested, say, in bounding $\|f +
h\|_{r}$ in terms of $\|f\|_s$ and $\|h\|_s$. Thus, there are
several possibilities to study the behavior of $\|f+h\|_{p}$  as $p$
changes.
\end{remark}

\end{document}